# Some remarks and conjectures about Hankel determinants of polynomials which are related to Motzkin paths

Johann Cigler

johann.cigler@univie.ac.at

**Abstract**

This note collects some results and conjectures for the generating functions of the Hankel determinants of certain polynomials which are related to Motzkin paths.

**1.** For non-negative integers $n,k$ let $\mathbf{M}_{n,k}$ be the set of all Motzkin paths from $(0,0)$ to $(n,k)$, i.e. all lattice paths from $(0,0)$ to $(n,k)$ consisting of up-steps $U=(1,1)$, down-steps $D=(1,-1)$ and horizontal steps $H=(1,0)$, which never run below the $x-$axis. For each path $P$ we define the weight $w(P)$ as the product of the weights of its steps, where the horizontal steps $H$ have weight $t$ and the up-steps $U$ and down-steps $D$ have weight 1.

Let $M_{n,k}(t) = \sum_{P \in \mathbf{M}_{n,k}} w(P)$ be the weight of all Motzkin paths from $(0,0)$ to $(n,k)$.

These weights satisfy

$$M_{n,k}(t) = M_{n-1,k-1}(t) + tM_{n-1,k}(t) + M_{n-1,k+1}(t) \tag{1.1}$$

with $M_{n,k}(t) = 0$ for $k < 0$ and $M_{0,k}(t) = [k=0]$.

For $k=0$ we get the Motzkin polynomials

$$M_n(t) = M_{n,0}(t) = \sum_{j=0}^{\lfloor \frac{n}{2} \rfloor} \binom{n}{2j} C_j t^{n-2j}. \tag{1.2}$$

The right-hand side follows because there are $\binom{n}{2j}$ ways to choose $2j$ positions for the up- and down-steps and $C_j = \frac{1}{j+1}\binom{2j}{j}$ ways to construct Dyck paths on them.

For $t=1$ we get the Motzkin numbers $M_n = M_n(1)$, $(M_n)_{n \geq 0} = (1,1,2,4,9,21,51,127,323,\ldots)$, (cf. [9], A001006), for $t=0$ the aerated Catalan numbers $(M_n(0))_{n \geq 0} = (1,0,1,0,2,0,5,0,\cdots)$, (cf. [9] A126120), and for $t=2$ the shifted Catalan numbers $(M_n(2))_{n \geq 0} = (1,2,5,14,\cdots)$, (cf. [9], A000108).

For positive $k$ we get the formula

$$M_{n,k}(t) = \sum_{j=0}^{\lfloor \frac{n-k}{2} \rfloor} a(k+2j,k) \binom{n}{2j+k} t^{n-2j-k}, \tag{1.3}$$

where $a(n,k)$ denotes the number of non-negative paths with up- and down-steps from $(0,0)$ to $(n,k)$. By the reflection principle we get $a(k+2j,k) = \binom{k+2j}{j} - \binom{k+2j}{j-1}$, i.e. the number of all paths $(0,0) \to (k+2j,k)$ minus the number of those paths which cross the $x-$axis. If we reflect the latter paths on the axis $x=-1$ after the first crossing we get a bijection with all paths $(0,0) \to (k+2j,-k-2)$, which have $j-1$ up-steps.

The first terms of the matrix $\left(M_{n,k}(t)\right)_{n,k \geq 0}$ are

$$\begin{pmatrix} 1 & 0 & 0 & 0 & 0 & 0 \\ t & 1 & 0 & 0 & 0 & 0 \\ 1+t^2 & 2t & 1 & 0 & 0 & 0 \\ t(3+t^2) & 2+3t^2 & 3t & 1 & 0 & 0 \\ 2+6t^2+t^4 & 4t(2+t^2) & 3(1+2t^2) & 4t & 1 & 0 \\ t(10+10t^2+t^4) & 5(1+4t^2+t^4) & 5t(3+2t^2) & 2(2+5t^2) & 5t & 1 \end{pmatrix}$$

For $t=0,1,2,3$ these triangles occur in OEIS[9] as A053121, A064189, A039598, A091965.

We want to obtain some information about the Hankel determinants

$$d_m^{(k)}(n,t) = \det\left(M_{m+i+j,k}(t)\right)_{i,j=0}^{n-1} \tag{1.4}$$

of the columns of this triangle.

For $k=0$ the determinants $d_m(n,t) = d_m^{(0)}(n,t)$ have been considered in [1], [2] and [4] with different methods. For $k>0$ and small $m$ the determinants $d_m^{(k)}(n,t)$ have been obtained and proved in [2]. We give an overview of these results and state conjectures for the remaining cases. In order to obtain nice results, we shall always set $d_m^{(k)}(0,t) = 1$.

The approach with orthogonal polynomials shows that the Motzkin polynomials $M_n(t)$ are the moments of the polynomials $p_n(x,t)$ which satisfy

$$p_n(x,t) = (x-t)p_{n-1}(x,t) - p_{n-2}(x,t) \tag{1.5}$$

with initial values $p_{-1}(x,t) = 0$ and $p_0(x,t) = 1$. Thus

$$p_n(x,t) = F_n(x-t), \tag{1.6}$$

where $F_n(x)$ are the Fibonacci polynomials defined by

$$F_n(x) = xF_{n-1}(x) - F_{n-2}(x) \tag{1.7}$$

with $F_0(x) = 1$ and $F_1(x) = x$.

This implies $d_0(n,t) = 1$ and

$$d_1(n,t) = (-1)^n p_n(0,t) = (-1)^n F_n(-t) = F_n(t). \tag{1.8}$$

The first terms are $d_1(0,t) = 1$, $d_1(1,t) = t$, $d_1(2,t) = t^2 - 1$, $d_1(3,t) = t^3 - 2t$, $d_1(4,t) = t^4 - 3t^2 + 1$.

The higher determinants can be computed using Dodgson's condensation.

The first terms of $d_2(n,t)$ are $1$, $t^2 + 1$, $t^4 - t^2 + 2$, $t^6 - 3t^4 + 3t^2 + 2$, $t^8 - 5t^6 + 8t^4 - 3t^2 + 3$.

There is no obvious regularity in the coefficients, but there are nice expressions in terms of Fibonacci polynomials.

In [1] it is shown that

$$d_2(n,t) = \sum_{j=0}^{n} F_j(t)^2 \tag{1.9}$$

and in [4] that

$$d_2(n,t) = \det \begin{pmatrix} F_n(t) & F_n'(t) \\ F_{n+1}(t) & F_{n+1}'(t) \end{pmatrix}. \tag{1.10}$$

For higher orders explicit expressions become more complicated, but there are nice generating functions.

Let

$$D_m(x,t) = \sum_{n \geq 0} d_m(n,t) x^n \tag{1.11}$$

denote the generating function of $d_m(n,t)$.

For small $m$ it can be proved that

$$D_0(x,t) = \frac{1}{1-x}, \tag{1.12}$$

$$D_1(x,t) = \frac{1}{1-tx+x^2}. \tag{1.13}$$

$$D_2(x,t) = \frac{1+x}{(1-x)^2 \left(1-\left(t^2-2\right)x+x^2\right)}, \tag{1.14}$$

$$D_3(x,t) = \frac{\left(1-x^2\right)\left(1+3tx+x^2\right)}{\left(1-(t^3-3t)x+x^2\right)\left(1-tx+x^2\right)^3}. \tag{1.15}$$

It turns out that there is a close relation with the Lucas polynomials $L_n(t)$. These are defined by

$$L_n(t) = tL_{n-1}(t) - L_{n-2}(t) \tag{1.16}$$

with initial values $L_0(t) = 2$ and $L_1(t) = t.$

Let

$$\alpha(t) = \frac{t + \sqrt{t^2 - 4}}{2},$$
$$\beta(t) = \frac{t - \sqrt{t^2 - 4}}{2} \tag{1.17}$$

be the roots of $x^2 - tx + 1 = 0.$

Binet's formulas show that

$$L_n(t) = \alpha(t)^n + \beta(t)^n,$$
$$F_n(t) = \frac{\alpha(t)^{n+1} - \beta(t)^{n+1}}{\alpha(t) - \beta(t)}. \tag{1.18}$$

Comparing the denominators $h_m(x,t)$ of $D_m(x,t)$ in the above examples we see that
$h_0(x,t) = 1 - x,$

$$h_1(x,t) = 1 - tx + x^2 = (\alpha(t) - x)(\beta(t) - x) = \left(1 - \frac{x}{\alpha(t)}\right)\left(1 - \frac{x}{\beta(t)}\right) = h_0\left(\frac{x}{\alpha(t)}, t\right) h_0\left(\frac{x}{\beta(t)}, t\right),$$

$$h_2(x,t) = h_1\left(\frac{x}{\alpha(t)}, t\right) h_1\left(\frac{x}{\beta(t)}, t\right) = (\alpha(t)\alpha(t) - x)(\alpha(t)\beta(t) - x)(\beta(t)\alpha(t) - x)(\beta(t)\beta(t) - x).$$

This leads to the conjecture that $D_m(x,t)$ can be written as a fraction with denominator
$h_m(x,t) = \prod_{i_1,\cdots,i_m} (\gamma_{i_1} \gamma_{i_2} \cdots \gamma_{i_m} - x),$ where the product is taken over all $(i_1, i_2, \cdots, i_m) \in \{0,1\}^m$ with
$\gamma_0 = \gamma_0(t) = \alpha(t), \quad \gamma_1 = \gamma_1(t) = \beta(t).$

I am indebted to Christian Krattenthaler for pointing out that a proof of this conjecture follows from his paper [8], Corollary 9, by setting there $s = t, \ t = 1$ and $x_i = 0.$

For $2j \neq m$ there are $\binom{m}{j}$ factors $(\alpha^{m-j}\beta^j - x) = (\alpha^{m-2j} - x)$ and $(\beta^{m-j}\alpha^j - x) = (\beta^{m-2j} - x).$

Therefore $h_m(x,t)$ has $\binom{m}{j}$ factors $(\alpha^{m-2j} - x)(\beta^{m-2j} - x) = x^2 - L_{m-2j}(t)x + 1.$

For even $m = 2\ell$ we get $\binom{m}{\ell}$ factors $1 - x.$

Setting

$$A_n(x,t) = x^2 - L_n(t)x + 1 = \left(x - \alpha(t)^n\right)\left(x - \beta(t)^n\right) \text{ for } n > 0, \quad (1.19)$$
$$A_0(x,t) = 1 - x,$$

we get $h_m(x,t) = \prod_{j=0}^{\lfloor\frac{m}{2}\rfloor} A_{m-2j}(x,t)^{\binom{m}{j}}$. This implies

**Theorem 1.1**

*For $m > 0$*

$$D_m(x,t) = \frac{R_m(x,t)}{\prod_{j=0}^{\lfloor\frac{m}{2}\rfloor} A_{m-2j}(x,t)^{\binom{m}{j}}} \quad (1.20)$$

with $R_m(x,t) \in \mathbb{Z}[x,t]$.

For $t \in \{0,1,2\}$ some (conjectured) simplifications are possible due to the periodicity of the sequence $\left(L_k(t)\right)_{k\geq 0}$.

For example, for $k = 0$ and $t = 1$ we get

$$D_{2m}(x,1) = \frac{r(2m,x)}{\left(1-x^3\right)^{m^2}(1-x)} \quad (1.21)$$

where $r(2m,x)$ is a palindromic polynomial with positive coefficients of degree $\deg r(2m,x) = m(3m-2)$ and

$$D_{2m+1}(x,1) = \frac{\left(1+(-1)^m x\right)(1+x)^2 r(2m+1,x,1)}{\left(1+x^3\right)^{m^2+m+1}} \quad (1.22)$$

where $r(2m+1,x)$ is a palindromic polynomial with degree $3m^2 + m + 1$.

For example, $r(1,x) = r(2,x) = 1 + x$, $r(3,x) = (1-x)(1+x)^2(1+3x+x^2)$,

$r(4,x) = 1 + 8x + 9x^2 + 14x^3 + 32x^4 + 14x^5 + 9x^6 + 8x^7 + x^8$,

$r(5,x) = (1+x)^3\left(1 + 18x + 9x^2 - 115x^3 - 203x^4 + 132x^5 + 384x^6 + 132x^7 - 203x^8 - 115x^9 + 9x^{10} + 18x^{11} + x^{12}\right)$.

Let us now consider the determinants $\left(d_m^{(k)}(n,t)\right)_{n\geq 0}$ for $k > 0$ and their generating functions

$$D_m^{(k)}(x,t) = \sum_{n\geq 0} d_m^{(k)}(n,t)x^n. \quad (1.23)$$

In order to stress the analogy with the case $k = 0$ we consider generalized Fibonacci polynomials $F_n^{(k)}(x)$ defined by

$$F_n^{(k)}(x) = L_{k+1}(x)F_{n-1}^{(k)}(x) - F_{n-2}^{(k)}(x) \tag{1.24}$$

with initial values $F_{-1}^{(k)}(x) = 0$ and $F_0^{(k)}(x) = 1$.

For $k = 0$ these reduce to $F_n^{(0)}(x) = F_n(x)$.

Since $\alpha^{k+1}$ and $\beta^{k+1}$ are roots of $(x - \alpha^{k+1})(x - \beta^{k+1}) = x^2 - L_{k+1}(t)x + 1$ we get as in (1.18)

$$F_n^{(k)}(t) = \frac{\alpha^{(k+1)(n+1)} - \beta^{(k+1)(n+1)}}{\alpha^{k+1} - \beta^{k+1}}. \tag{1.25}$$

By [2], Theorem 1, we know that $d_0^{(k)}((k+1)n, t) = (-1)^{n\binom{k+1}{2}}$ and $d_0^{(k)}(n,t) = 0$ else
and by [2], Theorem 2, and (1.25)

$$d_1^{(k)}((k+1)n, t) = (-1)^{n\binom{k+1}{2}} F_n^{(k)}(t), \quad d_1^{(k)}((k+1)n + k, t) = (-1)^{\binom{k}{2}} F_n^{(k)}(t), \text{ and } d_1^{(k)}(n,t) = 0 \text{ else.}$$

Therefore we have

$$D_0^{(k)}(x,t) = \sum_{n \geq 0} (-1)^{n\binom{k+1}{2}} x^{(k+1)n} = \frac{1}{1 - (-1)^{\binom{k+1}{2}} x^{k+1}} \tag{1.26}$$

and

$$D_1^{(k)}(x,t) = \frac{1 + (-1)^{\binom{k}{2}} x^k}{1 - (-1)^{\binom{k+1}{2}} L_{k+1}(t) x^{k+1} + x^{2(k+1)}}. \tag{1.27}$$

Let $h_m^{(k)}(x,t)$ be the denominator of $D_m^{(k)}(x,t)$. Then we have $h_0^{(k)}(x,t) = 1 - (-1)^{\binom{k+1}{2}} x^{k+1}$ and

$$h_1^{(k)}(x,t) = 1 - (-1)^{\binom{k+1}{2}} L_{k+1}(t) x^{k+1} + x^{2(k+1)} = \left(\alpha(t)^{k+1} - (-1)^{\binom{k+1}{2}} x^{k+1}\right)\left(\beta(t)^{k+1} - (-1)^{\binom{k+1}{2}} x^{k+1}\right)$$

$$= h_0^{(k)}\left(\frac{x}{\alpha(t)}, t\right) h_0^{(k)}\left(\frac{x}{\beta(t)}, t\right).$$

Further computations suggest that $D_m^{(k)}(x,t)$ can always be written as a fraction with denominator $h_m^{(k)}(x,t) = \prod_{i_1, \ldots, i_m} \left((\gamma_{i_1} \gamma_{i_2} \cdots \gamma_{i_m})^{k+1} - (-1)^{\binom{k+1}{2}} x^{k+1}\right).$

Setting

$$A_{k,0}(x,t) = 1 + (-1)^{\binom{k-1}{2}} x^{k+1},$$

$$A_{k,n}(x,t) = 1 + (-1)^{\binom{k-1}{2}} L_n(t) x^{k+1} + x^{2(k+1)}$$

(1.28)

we get

**Conjecture 1.2**

$$D_m^{(k)}(x,t) = \frac{R_m^{(k)}(x,t)}{\prod_{j=0}^{\lfloor \frac{m}{2} \rfloor} A_{k,(k+1)(m-2j)}^{\binom{m}{j}}(x,t)}$$

(1.29)

with $R_m^{(k)}(x,t) \in \mathbb{Z}[t,x]$.

A closer look at $D_m^{(k)}(x,t)$ suggests that some factors can be cancelled. This leads to

**Conjecture 1.3**

$$D_m^{(k)}(x,t) = \frac{r_m^{(k)}(x,t)}{\prod_{j=0}^{\lfloor \frac{m}{2} \rfloor} A_{k,(k+1)(m-2j)}^{1+j(m-j)}(x,t)}$$

(1.30)

where $r_m^{(k)}(x,t)$ is a polynomial of degree $\deg_x r_m^{(k)}(x,t) = \binom{m+1}{3} + k\left(\binom{m}{1} + \binom{m}{2} + \binom{m}{3}\right)$.

For $t = 2$ the matrix $\left(M_{n,k}(2)\right)_{n,k \geq 0}$ is the Catalan triangle [9], A039598, whose first terms are

$$\begin{pmatrix} 1 & 0 & 0 & 0 & 0 & 0 & 0 \\ 2 & 1 & 0 & 0 & 0 & 0 & 0 \\ 5 & 4 & 1 & 0 & 0 & 0 & 0 \\ 14 & 14 & 6 & 1 & 0 & 0 & 0 \\ 42 & 48 & 27 & 8 & 1 & 0 & 0 \\ 132 & 165 & 110 & 44 & 10 & 1 & 0 \\ 429 & 572 & 429 & 208 & 65 & 12 & 1 \end{pmatrix}.$$

Here we guess that

$$D_m^{(k)}(x,2) = \frac{A_m^{(k)}(x,2)}{\left(1 - (-1)^{\binom{k+1}{2}} x^{k+1}\right)^{1+\binom{m+1}{2}}}$$

(1.31)

where $A_m^{(k)}(x,2)$ is a polynomial of degree $\deg A_m^{(k)}(x,2) = \dfrac{m((k+1)m+k-1)}{2}$ which satisfies

$$x^{\frac{m((k+1)m+k-1)}{2}} A_m^{(k)}\left(\frac{1}{x},2\right) = \varepsilon(k,m) A_m^{(k)}(x,2) \tag{1.32}$$

with $\varepsilon(k,m)=1$ if $k \equiv 0 \pmod 4$, $\varepsilon(k,m)=(-1)^m$ if $k \equiv -1 \pmod 4$, $\varepsilon(k,m)=(-1)^{\binom{m}{2}}$ if $k \equiv 1 \pmod 4$, and $\varepsilon(k,m)=(-1)^{\binom{m+1}{2}}$ if $k \equiv 2 \pmod 4$.

For $k=0$ formula (1.31) is a known result. The Hankel determinants of the shifted Catalan numbers satisfy $\det\left(C_{m+i+j}\right)_{i,j=0}^{n-1} = \prod_{1 \leq i \leq j \leq m-1} \dfrac{2n+i+j}{i+j}$. I owe to Sam Hopkins [5] the observation that the right-hand side can be interpreted as the number of plane partitions of the form $(m-1, m-2, \cdots, 1)$ of non-negative integers $\leq n$. (A simple proof due to Christian Krattenthaler can be found in [3]). Formula (1.31) can be deduced from [10], Theorem 3.15.8.

The polynomials $A_m^{(0)}(x,2)$ are palindromic with positive coefficients of degree $\binom{m}{2}$. For example, $A_0^{(0)}(x,2)=1$, $A_1^{(0)}(x,2)=1$, $A_2^{(0)}(x,2)=1+x$, $A_3^{(0)}(x,2)=1+7x+7x^2+x^3$, $A_4^{(0)}(x,2)=1+31x+187x^2+330x^3+187x^4+31x^5+x^6$.

For $k>0$ we get for example

$$A_1^{(k)}(x,2) = 1 + (-1)^{\binom{k}{2}} x^k,$$

$$A_2^{(k)}(x,2) = 1 + (-1)^{\binom{k-1}{2}} x^{k-1} + (-1)^{\binom{k}{2}}(k+1)^2 x^k \left(1 + (-1)^{\binom{k+1}{2}} x^{k+1}\right) - x^{2k+2} \text{ for } k \geq 2, \text{ and}$$

$$A_2^{(1)}(x,2) = (1-x^2)(1+4x+x^2).$$

Apparently there are also analogs of (1.9) and (1.10):

**Conjecture 1.4**

$$(-1)^{\binom{k+1}{2}n} d_2^{(k)}\left((k+1)n, t\right) = F_n^{(k)}(t)^2,$$

$$(-1)^{\binom{k+1}{2}n + \binom{k-1}{2}} d_2^{(k)}\left((k+1)n+k-1, t\right) = F_n^{(k)}(t)^2, \tag{1.33}$$

$$(-1)^{\binom{k+1}{2}n + \binom{k}{2}} d_2^{(k)}\left((k+1)n+k, t\right) = L_{k+1}(t)' \sum_{j=0}^{n} F_j^{(k)}(t)^2 = \det\begin{pmatrix} F_n^{(k)}(t) & F_n^{(k)}(t)' \\ F_{n+1}^{(k)}(t) & F_{n+1}^{(k)}(t)' \end{pmatrix},$$

$d_2^{(k)}(n,t)=0$ else.

**2.** Now let us consider a slight generalization by changing the weight of the horizontal steps $H$ on height $0$ to $s$ instead of $t$. Let $M_{n,k}(t,s)$ denote these weights of the Motzkin paths.

Here we get $d_0(n,t,s) = 1$ and

$$\sum_{n \geq 0} d_1(n,t,s)x^n = \frac{1+(s-t)x}{1-tx+x^2} \tag{2.1}$$

$$\sum_{n \geq 0} d_2(n,t,s)x^n = \frac{1+\left(1+s^2-t^2\right)x+(s-t)^2 x^2}{(1-x)^2 \left(1+\left(2-t^2\right)x+x^2\right)}. \tag{2.2}$$

**Conjecture 2.1**

$$\sum_{n \geq 0} d_m(n,t,s)x^n = \frac{R_m(x,t,s)}{\prod_{j=0}^{\left\lfloor \frac{m}{2} \right\rfloor} A_{0,m-2j}^{1+j(m-j)}(x,t)} \tag{2.3}$$

where $R_m(x,t,s)$ is a polynomial in $x, s, t$ with integer coefficients with

$$\deg_x R_m(x,t,s) = \binom{m+1}{3}+1.$$

**Conjecture 2.2**

For all $s$ the sequences $\left(d_m^{(k)}(n,t,0)\right)_{n \geq 0}$ and $\left(d_m^{(k)}(n,t+s,s)\right)_{n \geq 0}$ satisfy the same recurrence of order $2^{k+m}$ for all $k, m > 0$. In the special case $m = 0$ we even get

$$d_0^{(k)}(n,t+s,s) = d_0^{(k)}(n,t,0). \tag{2.4}$$

Let us mention some explicit formulas for some small $m$ and $k$:

$$\sum_{n \geq 0} d_0^{(1)}(n,t,0)x^n = \frac{1-tx}{1-tx+x^2} \tag{2.5}$$

$$\sum_{n \geq 0} d_0^{(2)}(n,t,0)x^n = \frac{1+x+t^2 x^2}{1+x+t^2 x^2 + x^3 + x^4} \tag{2.6}$$

$$\sum_{n \geq 0} d_0^{(3)}(n,t,0)x^n = \frac{1+tx+6\binom{t+1}{3}x^3 + \left(t^4+t^2-1\right)x^4 + t^3 x^5}{1+tx+6\binom{t+1}{3}x^3 + \left(t^2-1\right)\left(t^2+2\right)x^4 + 6\binom{t+1}{3}x^5 + tx^7 + x^8}. \tag{2.7}$$

The denominator for $k = 4$ is

$$1-x+t^2 x^3 - t^2 \left(2t^2-1\right)x^4 + \left(t^2-1\right)\left(t^2+3\right)x^5 - \left(2t^2-3\right)x^6 - t^2\left(t^4-t^2+1\right)x^7 + \left(t^2-1\right)t^2 \left(t^4+t^2+2\right)x^8$$
$$-t^2\left(t^4-t^2+1\right)x^9 - \left(2t^2-3\right)x^{10} + \left(t^2-1\right)\left(t^2+3\right)x^{11} - t^2\left(2t^2-1\right)x^{12} + t^2 x^{13} - x^{15} + x^{16}$$

$$\sum_{n\geq 0} d_1^{(1)}(n,t,s)x^n = \frac{1+(1+t(s-t))x+(s-t)^2 x^2}{1+(s-t)x+(L_2(t)+(s-t)^2)x^2+(s-t)x^3+x^4} \qquad (2.8)$$

$$\sum_{n\geq 0} d_1^{(2)}(n,t,0)x^n$$

$$= \frac{1+tx+t^2(t^2-2)x^2+t(2t^2-3)x^3+(t^4+t^2-1)x^4+t^3 x^5}{1+tx+t^2(t^2-2)x^2+t(2t^2-3)x^3+t^2(2t^2-3)x^4++t(2t^2-3)x^5+t^2(t^2-2)x^6+tx^7+x^8}.$$
(2.9)

For $(t,s)=(2,1)$ we get another Catalan triangle (cf. [9], A039599), whose first terms are

$$\begin{pmatrix} 1 & 0 & 0 & 0 & 0 & 0 & 0 \\ 1 & 1 & 0 & 0 & 0 & 0 & 0 \\ 2 & 3 & 1 & 0 & 0 & 0 & 0 \\ 5 & 9 & 5 & 1 & 0 & 0 & 0 \\ 14 & 28 & 20 & 7 & 1 & 0 & 0 \\ 42 & 90 & 75 & 35 & 9 & 1 & 0 \\ 132 & 297 & 275 & 154 & 54 & 11 & 1 \end{pmatrix}.$$

In this case we get more information about the generating functions:

$$\sum_{n\geq 0} d_m^{(k)}(n,2,1)x^n = \frac{a_m^{(k)}(x)}{\left(1-(-1)^k x^{2k+1}\right)^{\binom{m}{2}+1}}, \qquad (2.10)$$

with $\deg a_m^{(k)}(x) = (m(m-1)+1)k + \binom{m-1}{2}$ and

$$\mathrm{sgn}\left(\frac{a_m^{(k)}\left(\frac{1}{x}\right)}{a_m^{(k)}(x)}\right) = (-1)^{\binom{k}{2}} \text{ if } k \equiv 1,2 (\mathrm{mod}\, 4) \text{ and } \mathrm{sgn}\left(\frac{a_m^{(k)}\left(\frac{1}{x}\right)}{a_m^{(k)}(x)}\right) = (-1)^{\binom{k+1}{2}} \text{ if } k \equiv 3,0 (\mathrm{mod}\, 4).$$

## 3. Addendum

**3.1.** As another generalization we define for $P \in \mathbf{M}_{n,k}$ a weight $w_{t,s}(P)$ as the product of the weights of its steps, where $w_{t,s}(H)=t$, $w_{t,s}(D)=s$ and $w_{t,s}(U)=1$.

Let $M_{n,k,s}(t) = \sum_{P\in \mathbf{M}_{n,k}} w_{t,s}(P)$ be the weight of all paths from $(0,0)$ to $(n,k)$.

These weights satisfy

$$M_{n,k,s}(t) = M_{n-1,k-1,s}(t) + tM_{n-1,k,s}(t) + sM_{n-1,k+1,s}(t) \qquad (3.1)$$

with $M_{n,k,s}(t) = 0$ for $k < 0$ and $M_{0,k,s}(t) = [k=0]$.

The same argument as above gives

$$M_{n,k,s}(t) = \sum_{j=0}^{\lfloor \frac{n-k}{2} \rfloor} \binom{n}{2j+k} \binom{k+2j}{j} \frac{k+1}{k+j+1} t^{n-2j-k} s^j. \tag{3.2}$$

We want to give a generalization of Conjecture 1.3 to this case. In the same way as above, we see that

$$\det \left( M_{i+j,0,s}(t) \right)_{i,j=0}^{n-1} = s^{\binom{n}{2}}, \tag{3.3}$$

$$\det \left( M_{1+i+j,0,s}(t) \right)_{i,j=0}^{n-1} = s^{\binom{n}{2}} \sqrt{s}^n F_n\left( \frac{t}{\sqrt{s}} \right). \tag{3.4}$$

Therefore, we will consider instead of $\det \left( M_{m+i+j,k,s}(t) \right)_{i,j=0}^{n-1}$ the ratios

$$S_m^{(k)}(n,t,s) = \frac{\det \left( M_{m+i+j,k,s}(t) \right)_{i,j=0}^{n-1}}{s^{\binom{n}{2}}} \tag{3.5}$$

with $S_m^{(k)}(0,t,s) = 1$. From (3.4) we get

$$\sum_{n \geq 0} S_1^{(0)}(n,t,s) = \frac{1}{1 - tx + sx^2}, \tag{3.6}$$

where $1 - tx + sx^2 = (1 - \alpha(t,s)x)(1 - \beta(t,s)x)$ with

$$\begin{aligned} \alpha(t,s) &= \frac{t + \sqrt{t^2 - 4s}}{2}, \\ \beta(t,s) &= \frac{t - \sqrt{t^2 - 4s}}{2}. \end{aligned} \tag{3.7}$$

Let

$$L_n(t,s) = \alpha(t,s)^n + \beta(t,s)^n \tag{3.8}$$

and

$$A_{k,n}(x,t,s) = \left( 1 - \frac{\alpha(t,s)^n}{(-s)^{\binom{k+1}{2}}} x^{k+1} \right)\left( 1 - \frac{\beta(t,s)^n}{(-s)^{\binom{k+1}{2}}} x^{k+1} \right) = 1 - \frac{L_n(t,s)}{(-s)^{\binom{k+1}{2}}} x^{k+1} + \frac{s^n}{s^{k(k+1)}} x^{2(k+1)},$$

$$A_{k,0}(x,t,s) = 1 - \frac{x^{k+1}}{(-s)^{\binom{k+1}{2}}}. \tag{3.9}$$

Then computer experiments suggest

**Conjecture 3.1**

$$\sum_{n \geq 0} S_m^{(k)}(n,t,s)x^n = \frac{U_m^{(k)}(x,t,s)}{\prod_{j=0}^{\lfloor m/2 \rfloor} A_{k,(k+1)(m-2j)}^{1+j(m-j)}(s^j x, t, s)} \tag{3.10}$$

where $U_m^{(k)}(x,t,s)$ has degree $\deg_x U_m^{(k)}(x,t,s) = \binom{m+1}{3} + k\left(\binom{m}{1} + \binom{m}{2} + \binom{m}{3}\right).$

**3.2.** Finally let us consider the set $\mathbf{T}_n$ of ALL lattice paths from $(0,0)$ to $(n,0)$ consisting of up-steps $U = (1,1)$, down-steps $D = (1,-1)$ and horizontal steps $H = (1,0)$ with the weight $w_{t,s}$. Then we get

$$w_{t,s}(\mathbf{T}_n) = T_n(t,s) = \sum_{j=0}^{\lfloor n/2 \rfloor} \binom{n}{2j}\binom{2j}{j} t^{n-2j} s^j = [x^n](1+tx+sx^2)^n. \tag{3.11}$$

The right-hand identity follows from

$$(1+tx+sx^2)^n = \left(1+tx\left(1+\frac{sx}{t}\right)\right)^n = \sum_k \binom{n}{k}(tx)^k\left(1+\frac{sx}{t}\right)^k = \sum_{k,j}\binom{n}{k}\binom{k}{j}t^{k-j}s^j x^{k+j}$$

$$= \sum_n x^n \sum_j \frac{n!}{(j)!j!(n-2j)!} t^{n-2j} s^j = \sum_n x^n \sum_{j=0}^{\lfloor n/2 \rfloor} \binom{n}{2j}\binom{2j}{j} t^{n-2j} s^j.$$

The numbers $T_n = |\mathbf{T}_n|$ of these paths are therefore the central trinomial coefficients

$$T_n = [x^n](1+x+x^2)^n = \sum_{j=0}^{\lfloor n/2 \rfloor} \binom{n}{2j}\binom{2j}{j}. \tag{3.12}$$

Let now $\mathbf{d}_m(n,t,s) = \det\left(T_{m+i+j}(t,s)\right)_{i,j=0}^{n-1}$ with $\mathbf{d}_m(0,t,s) = 1$ for $m > 0$ and $\mathbf{d}_0(0,t,s) = \frac{1}{2}$.

Then $\mathbf{d}_0(n,t,s) = 2^{n-1} s^{\binom{n}{2}}$ for $n \in \mathbb{N}$.

For

$$\mathbf{D}_m(x,t,s) = \sum_{n \geq 0} \frac{\mathbf{d}_m(n,t,s)}{2\mathbf{d}_0(n,t,s)} x^n = \sum_{n \geq 0} \frac{\mathbf{d}_m(n,t,s)}{2^n s^{\binom{n}{2}}} x^n \tag{3.13}$$

we get

**Conjecture 3.2**

$$\mathbf{D}_m(x,t,s) = \frac{\mathbf{U}_m(x,t,s)}{\prod_{j=0}^{\lfloor \frac{m}{2} \rfloor} A_{0,m-2j}^{1+j(m-j)}(s^j x, t, s)} \tag{3.14}$$

with $\deg_x \mathbf{U}_m(x,t,s) = \binom{m}{3} + \binom{m-1}{2} + m + 1 = \frac{m^3 - m + 12}{6}$. for $m > 0$.

For $T_n(2,1) = \binom{2n}{n}$ this simplifies to

**Conjecture 3.3**

$$\mathbf{D}_m(x,2,1) = \frac{\mathbf{V}_m(x,2,1)}{(1-x)^{\binom{m}{2}+1}}, \tag{3.15}$$

where $\mathbf{V}_m(x,2,1)$ is a symmetric polynomial with positive coefficients and
$\deg \mathbf{V}_m(x,2,1) = \binom{m-1}{2} + 1$ for $m > 0$.

The first numerators are $\mathbf{V}_1(x,2,1) = 1$, $\mathbf{V}_2(x,2,1) = 1+x$, $\mathbf{V}_3(x,2,1) = 1+6x+x^2$, $\mathbf{V}_4(x,2,1) = 1+28x+70x^2+28x^3+x^4$.

The computations have been made with Mathematica and the Mathematica packages Guess by Manuel Kauers [6] and RATE by Christian Krattenthaler [7].